\documentclass[aap]{imsart}
\RequirePackage[OT1]{fontenc}
\RequirePackage{amsthm,amsmath}
\RequirePackage[numbers]{natbib}
\RequirePackage[colorlinks,citecolor=blue,urlcolor=blue]{hyperref}
\usepackage{mathrsfs,url,color}
\usepackage{latexsym}
\usepackage{amsfonts}
\usepackage{changes}
\usepackage{amsmath,amsfonts,amsthm}
\usepackage{mathrsfs}
\usepackage{amssymb, color}
\usepackage{pgf,tikz}
\usepackage{pgf,tikz}
\usetikzlibrary{decorations.pathreplacing}
\usetikzlibrary{calc}
\usepackage{pgf} 
\usepackage[english]{babel}
\usepackage[utf8]{inputenc}
\usepackage{color}
\usepackage{amssymb,color}

\startlocaldefs
\numberwithin{equation}{section}
\theoremstyle{plain}
\newtheorem{thm}{Theorem}[section]
\newtheorem{rem}[thm]{Remark}

\newtheorem{lem}[thm]{Lemma}

\newtheorem{obs}[thm]{facts}

\newcommand{\lemref}[1]{Lemma~{\rm \ref{#1}}}

\endlocaldefs

\begin{document}
	
	\begin{frontmatter}
		\title{Optimal Bounds in Normal Approximation for Many Interacting Worlds}
		\runtitle{Normal Approximation for Many Interacting Worlds}
		
		\begin{aug}
			\author{\fnms{Louis H. Y.} \snm{Chen}\thanksref{m1}\ead[label=e1]{matchyl@nus.edu.sg}}
			\and
			\author{\fnms{L\^{e} V\v{a}n} \snm{Th\`{a}nh}\thanksref{m2}\ead[label=e2]{levt@vinhuni.edu.vn}}

			\runauthor{Louis H. Y. Chen and L\^{e} V\v{a}n Th\`{a}nh}
			
			\affiliation{National University of Singapore\thanksmark{m1} and Vinh University\thanksmark{m2}}
			
			\address{Department of Mathematics\\
				National University of Singapore\\
				10 Lower Kent Ridge Road\\ 
				Singapore 119076}

			\address{Department of Mathematics\\
				Vinh University\\
				182 Le Duan, Vinh, Nghe An\\
				Vietnam}
				
		\end{aug}
		
		\begin{abstract}. In this paper, we use Stein's method to obtain optimal bounds, both in Kolmogorov and in Wasserstein distance, in the normal approximation for the empirical distribution of the ground state of a many-interacting-worlds harmonic oscillator proposed by Hall, Deckert, and Wiseman [{\it Phys. Rev. X.} (2014)]. Our bounds on the Wasserstein distance solve a conjecture of McKeague and Levin [{\it Ann. Appl. Probab.} (2016)].
		\end{abstract}
		\begin{keyword}[class=MSC]
			\kwd[Primary ]{60F05, 81Q65}
		\end{keyword}
		
		\begin{keyword}
			\kwd{Stein's method}
			\kwd{normal approximation}
			\kwd{many interacting worlds}
			\kwd{harmonic oscillator}
		\end{keyword}
		
	\end{frontmatter}
	
	\section{Introduction}\label{sec:int}

	In \cite{HDW14}, Hall, Deckert, and Wiseman proposed a many-interacting-worlds (MIW) theory for interpreting quantum mechanics. In this theory, quantum theory can be understood as the continuum limit of a deterministic theory in which there is a large, but finite, number of interacting classical ``worlds''. Here, a world means an entire universe with well-defined properties, determined by the classical configuration of its particles and fields.

	Hall, Deckert, and Wiseman \cite{HDW14} proposed a MIW harmonic oscillator model for $N$ one-dimensional worlds, where the Hamiltonian for the MIW harmonic oscillator is:
$$
H(\mathbf{x}, \mathbf{p}) =  E(\mathbf{p}) + V(\mathbf{x}) + U(\mathbf{x}).
$$
	Here $\mathbf{p}=(p_1,\ldots,p_N)$ are the momenta and $\mathbf{x} = (x_1, \cdots, x_N)$, $x_1 > x_2 > \cdots > x_N$, the locations of the $N$ worlds, each regarded as a particle with unit mass,
	$$E(\mathbf{p}) = \sum_{n=1}^N p_n^2/2$$
	is the kinetic energy,
	$$V(\mathbf{x}) = \sum_{n=1}^N x_n^2$$
	is the potential energy, and 
	$${\displaystyle U(\mathbf{x}) = \sum_{n=1}^N \left(\frac{1}{x_{n-1}-x_n} - \frac{1}{x_n - x_{n+1}}\right)^2}$$
	is the ``interworld'' potential, where $x_0 = \infty$ and $x_{N+1} = -\infty$, which is a discretization of Bohm's quantum potential (see Bohm \cite{Bohm1,Bohm2}).
	
	Hall, Deckert, and Wiseman \cite{HDW14} showed that in the ground state, where the Hamiltonian is minimized, all the momenta $p_n$ vanish and the locations of the $N$ particles satifisfy this recursion equation:
	\begin{equation}
		\label{MIW-recursion}
		x_{n+1} = x_n - \frac{1}{x_1 + \cdots + x_n}, ~ 1\le n\le N-1,
	\end{equation}
	suject to the constraints $x_1+\cdots+ x_N = 0$ and $x_1^2+\cdots x_N^2 = N - 1$. 
	
	Their numerical calculations suggest that the empirical distribution ${\mathbb{P}}_N$ of the locations converges to the standard Gaussian distribution $\gamma$ as $N \rightarrow \infty$, which agrees with ground state probability distribution of a quantum harmonic oscillator. Here the empirical distribution ${\mathbb{P}}_N$ is defined by
	\begin{equation}\label{empirical}
		{\mathbb{P}}_N(A) :=\dfrac{\#\{n:x_n\in A\}}{N} ~ \text{for}~ A \in \mathcal{B}(\mathbb{R}).
	\end{equation}
	
	In \cite{ML16}, McKeague and Levin proved that the recursion equation \eqref{MIW-recursion} has a unique solution if the solution is
	monotonic, zero-median, namely $x_m = 0$ for $N$ odd and $x_m = - x_{m+1}$ for $N$ even, where $m=(N+1)/2$ if $N$ is odd, and $m=N/2$ if $N$ is even. They proved that ${\mathbb{P}}_N$ converges to $\gamma$ as $N \rightarrow \infty$.

The Wasserstein distance between ${\mathbb{P}}_N$ and $\gamma$, denoted by $d_{\mathrm{W}}({\mathbb{P}}_N,\gamma)$, is defined by
	\[d_{\mathrm{W}}({\mathbb{P}}_N, \gamma):= \sup_{|h(x)-h(y)\le|x-y|}\left|\int_\mathbb{R}hd{\mathbb{P}}_N -\int_\mathbb{R}hd\gamma\right|.\]
	Using Stein's method and zero-bias coupling, McKeague and Levin \cite{ML16} further proved that 
	\[d_{\mathrm{W}}({\mathbb{P}}_N,\gamma) \le \frac{4}{\sqrt{\log N}},\]
	and conjectured that the correct order of the bound on $d_{\mathrm{W}}({\mathbb{P}}_N,\gamma)$ should be $\sqrt{\log N}/N$.

The Kolmogorov distance between ${\mathbb{P}}_N$ and $\gamma$, denoted by $d_{\mathrm{K}}(\mathbb{P}_N, \gamma)$, is defined by
	\[d_{\mathrm{K}}({\mathbb{P}}_N, \gamma):= \sup_{z\in \mathbb{R}}\left|\mathbb{P}_N((-\infty,z]) - \gamma((-\infty,z])\right|.\]

In this paper, we prove that the upper and lower bounds on the Wasserstein distance between $\mathbb{P}_N$ and $\gamma$  are both of the order $\sqrt{\log N}/N$, thereby proving the conjecture of McKeague and Levin \cite{ML16} while at the same time showing that $\sqrt{\log N}/N$ is optimal. We also prove that the upper and lower bounds on the Kolmogorov distance between $\mathbb{P}_N$ and $\gamma$ are both of the order $1/N$, showing that $1/N$ is also optimal. This is a surprising outcome as optimal bounds on the Kolmogorov distance are usually of no smaller order than those on the Wasserstein distance for a particular problem.  

Our proof of the upper bound on the Wasserstein distance turned out to be an easy consequence of Theorem 1.1 of Goldstein \cite{Goldstein07} using the zero-bias coupling of McKeague and Levin \cite{ML16} and an upper bound on $x_1$. Our approach is also the same as that of the proof of the upper bound on the Wasserstein distance for two-sided Maxwell approximation in McKeague, Pek\"{o}z and Swan \cite[p. 109]{MPS19}.
	
\begin{rem} \rm Stein's method is applicable to Wasserstein distance of order $p > 1$ (see, for example, \cite{Bonis,CFP,Fang,LNP} for recent results in this direction). We have not succeeded in obtaining optimal bounds on the Wasserstein distance of order $p > 1$ between $\mathbb{P}_N$ and $\gamma$. It remains an open problem for future research.
\end{rem}

	Throughout this paper,
	$\gamma$ denotes the standard normal distribution.
	Let $N\ge 2$ be an integer number, we will denote 
	$m=(N+1)/2$ if $N$ is odd, and $=N/2$ if $N$ is even.
	For a positive number $x$, $\lfloor x \rfloor$ denotes the greatest integer number which is less than or equal to $x$, and $\log x$ denotes the natural (base $e= 2.7182\ldots$) logarithm of $x$.
	The symbol $C$ denotes a positive constant which does not depend on $N$, and its value
	may be different at each appearance. 
	We also denote by
	$\{x_1,\dots,x_N\}$ the unique zero-median and strictly decreasing
	solution 
	to the recursion equation \eqref{MIW-recursion} and by ${\mathbb{P}}_N$ the empirical distribution as in \eqref{empirical}, and
	let
	\[S_n=\sum_{i=1}^n x_i,\ 1\le n\le N.\]

	The rest of the paper is organized as follows.
	In Section 2, we reproduce the construction of the zero-bias coupling of 
	McKeague and Levin \cite{ML16} and state the main result. 
	Section 3 focuses on
	the proof of the bounds 
	on the Kolmogorov distance.
	The bounds on the Wasserstein distance are presented in Section 4. 
	Finally, some technical results are proved in the Appendix.

	\section{Main results}\label{sec:main}
	
	Before proving the main result, we will need some preliminary lemmas. The first lemma establishes the existence
	of a solution to the recursion relation.  This lemma is proved in
	McKeague and Levin \cite{ML16}. 
	
\begin{lem}\label{ML1}
		Every zero-median solution $\{x_1,\dots,x_N\}$ of \eqref{MIW-recursion} satisfies the following properties: 
		
		(P1) Zero-mean: $x_1 + \dots + x_N = 0$.
		
		(P2) Variance-bound: $x_{1}^{2} + \dots+x_{N}^{2}= N-1.$
		
		(P3) Symmetry: $x_n =  -x_{N+1-n}$ for $n = 1,\dots,N.$
		
		Further, there is a unique solution $\{x_1,\dots,x_N\}$ of \eqref{MIW-recursion}
		such that (P1) holds and
		
		(P4) Strictly decreasing: $x_1 > \dots > x_N$.
		
		This solution has the zero-median property, and thus also satisfies (P2)
		and (P3).
	\end{lem}
	
	It was shown by Goldstein and Reinert \cite{GoldsteinReinert97} that
	for any mean zero random variable $W$ with finite variance
	$\sigma^2$, there exists a random variable $W^{*}$ which satisfies
	\begin{equation}\label{zero}
		\mathbb{E}Wf(W)=\sigma^2\mathbb{E}f^{'}(W^{*})
	\end{equation}
	for all absolutely
	continuous $f$ with $\mathbb{E}|Wf(W)|<\infty$. We say such a $W^{*}$
	has the $W$-zero-biased distribution. The following result is due to
	Goldstein and Reinert \cite{GoldsteinReinert97} (see also in
	Chen, Goldstein and Shao \cite[Proposition 2.1]{CGS}).
	\begin{lem}\label{GR}
		Let $W$ be a random variable with mean zero and finite positive
		variance $\sigma^2$, and let $W^{*}$ have the $W$-zero-biased
		distribution. Then the distribution of $W^{*}$ is absolutely
		continuous with density given by
		$$p^*(x)=\mathbb{E}\left(W{\mathbf{1}}(W>x)\right)/\sigma^2=-\mathbb{E}\left(W{\mathbf{1}}(W\le x)\right)/\sigma^2.$$
	\end{lem}
	We recall that ${\mathbb{P}}_N$ denotes the empirical probability measure as in \eqref{empirical}.
	McKeague and Levin \cite[page 9]{ML16} constructed a zero bias coupling $(W,W^*)$,
	where $W$ has the probability distribution ${\mathbb{P}}_N$. For completeness, we describe the construction in this paper.
	From \lemref{ML1}, we have $\operatorname{Var}(W)=(N-1)/N$. By \lemref{GR}, the density of $W^*$ is given by
	$$p^*(x)=\dfrac{S_n}{N-1}=\dfrac{1}{(N-1)(x_n-x_{n+1})}\ \mbox{ if } x_{n+1}\le x< x_n,\ 1 \le n \le N-1.$$
	This implies that $W^*$ is uniformly distributed on each interval $[x_{n+1},x_n)$ with mass $1/(N-1)$, $1 \le n \le N-1$. Let $x_{n+1} < y_n < x_n$ such that for $1 \le n\le N-1$, the area under $p^*$ on the interval $[y_n,x_n)$ is $L_n = (N-n)/(N(N-1))$, and for $2 \le n \le N$, the area under $p^*$ on the interval $[x_n,y_{n-1})$ is $R_n = (n-1)/(N(N-1))$. Then the area under $p^*$ on $[y_1,x_1)$ and on $[x_N,y_{N-1})$ is $1/N$, and on $[y_n,y_{n-1})$ is $L_n+R_n = 1/N$ for $2 \le n \le N-1.$  See Figure 1 below.

\vskip.3in
\begin{figure}[h!]
\def\yscale{0.1}
\def\bluerecright#1#2#3#4{
{\color{blue}
\def\i{#1}
\pgfmathparse{int(\i+1)}
\edef\iplusone{\pgfmathresult}
\pgfmathparse{#2*\xns[\i]+#3*\xns[\iplusone]}
\edef\x{\pgfmathresult}
\pgfmathparse{divide(\yscale,\xns[\i]-\xns[\iplusone])}
\edef\y{\pgfmathresult}
\filldraw (\x,0) rectangle (\xns[\i],\y)
 [densely dotted,line width=0.7pt,fill=#4];}
}
\def\bluerecleft#1#2#3#4{
{\color{blue}
\def\i{#1}
\pgfmathparse{int(\i+1)}
\edef\iplusone{\pgfmathresult}
\pgfmathparse{#2*\xns[\i]+#3*\xns[\iplusone]}
\edef\x{\pgfmathresult}
\pgfmathparse{divide(\yscale,\xns[\i]-\xns[\iplusone])}
\edef\y{\pgfmathresult}
\filldraw (\x,0) rectangle (\xns[\iplusone],\y)
 [densely dotted,line width=0.7pt,fill=#4];}
}
\def\bluelabel#1#2#3#4{
{\color{blue}
\def\i{#1}
\pgfmathparse{int(\i+1)}
\edef\iplusone{\pgfmathresult}
\pgfmathparse{#3*\xns[\i]+#4*\xns[\iplusone]}
\edef\x{\pgfmathresult}
\draw (\x,0) -- (\x,-0.2) [densely dotted,line width=0.7pt];
\draw (\x,-0.2) node[below=-1pt] {\footnotesize$y_{#2}$};
}}

\vspace{-5ex}

\begin{center}
\begin{tikzpicture}[scale=3]
\def\xns{{1.645,  1.118,  0.794,  0.538,  0.313,  0.103,
         -0.103, -0.313, -0.538, -0.794, -1.118, -1.645}}
         
\bluerecright{0}{0.25}{0.75}{blue!15}
\bluelabel{0}{1}{0.25}{0.75}
\bluerecleft{3}{0.5}{0.5}{blue!15}
\bluerecright{4}{0.5}{0.5}{blue!15}
\bluelabel{4}{n}{0.5}{0.5}
\bluelabel{3}{n-1}{0.5}{0.5}
\bluerecleft{10}{0.72}{0.28}{blue!15}
\bluelabel{10}{N-1}{0.72}{0.28}
\draw[-latex,line width=1.5pt] (-1.8,0) -- (1.8,0);
\foreach \i in {0,...,11} {
  \pgfmathparse{\xns[\i]}
  \node (A\i) at (\pgfmathresult,0) {};
  \draw[line width=1pt] ($(A\i)$) -- ($(A\i)-(0,0.03)$);
}
\foreach \i in {0,...,10} {
  \pgfmathparse{int(\i+1)}
  \edef\iplusone{\pgfmathresult}
  \pgfmathparse{\yscale/(\xns[\i]-\xns[\iplusone])}
  \draw ($(A\i)$) rectangle ($(A\iplusone)+(0,\pgfmathresult)$)
   [line width=1pt];
}

{\footnotesize
\draw  (A0) node[below=5pt] {$x_{1}$};
\draw  (A1) node[below=5pt] {$x_{2}$};
\draw (A10) node[below=5pt] {$x_{N-1}$};
\draw (A11) node[below=5pt] {$x_{N}$};
\draw  (A3) node[below=5pt] {$\,\,x_{n-1}$};
\draw  (A4) node[below=5pt] {$x_{n}$};
\draw  (A5) node[below=5pt] {$x_{n+1}\,\,$};
}
  
\draw[-latex,blue](0.25,0.38) .. controls +(90:2mm) and +(-20:2mm) .. (-0.1,0.7) node[left=0pt] 
{\color{blue}\footnotesize $\displaystyle L_n=\frac{N-n}{N(N-1)}$} ;

\draw[-latex,blue](0.36,0.38) .. controls +(90:2mm) and +(180:2mm) .. (0.75,0.7) node[right=0pt] 
{\color{blue}\footnotesize $\displaystyle
 R_n=\frac{n-1}{N(N-1)}$} ;

\draw[-latex,blue] (1.45,0.1) .. controls +(90:1mm) and +(210:1mm) .. (1.6,0.3) node [right=0pt] 
 {\footnotesize $L_1=\displaystyle\frac{1}{N}$};

\draw[-latex,blue] (-1.45,0.1) .. controls +(90:1mm) and +(-30:1mm) .. (-1.6,0.3) node [left=0pt] 
  {\footnotesize $R_{N}=\displaystyle\frac{1}{N}$};

\end{tikzpicture}

\end{center}
\caption{Density $p^*$ of $W^*$}
\label{fig:coupling}
\end{figure}
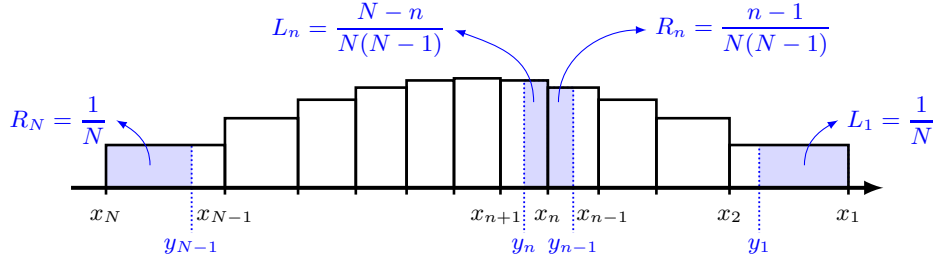	
	Now we couple $W$ to $W^*$ by defining $W$ and $W^*$ on the same probability space with  $\Omega=[x_N,x_1)$ as the sample space and $p^*$ as the probability measure as follows: $W^*(w)=\omega$ for all $\omega\in\Omega$, and
	\begin{equation*}
		W(w)=
		\begin{cases}
			x_1& \mbox{ if }\omega \in [y_1,x_1),\\
			x_n& \mbox{ if }\omega \in [y_{n},y_{n-1}), \mbox{  } 2\le n\le N-1,\\
			x_N& \mbox{ if }\omega \in [x_N,y_{N-1}).
		\end{cases}
	\end{equation*}
	
	The following theorem is the main result of the paper.

	\begin{thm}\label{main1}
		Let $\{x_n,1\le n\le N\}$ be the unique monotonic zero-mean solution of the recursion relation \eqref{MIW-recursion}
		and ${\mathbb{P}}_N$ the empirical distribution. Then
		\begin{equation}
			\label{BE1}
			\dfrac{1}{2N}\le d_K({\mathbb{P}}_N,\gamma)\le \dfrac{55}{N},
		\end{equation}
		and
		\begin{equation}
			\label{BE2}
			\dfrac{\sqrt{\log (N/2)}}{2N}-\dfrac{C}{N} \le d_W({\mathbb{P}}_N,\gamma)\le \dfrac{16\sqrt{\log N}}{N}.
		\end{equation}
	\end{thm}

	\section{Bounds on the Kolmogorov distance}\label{sec:Kolmogorov}

	In this section, we will prove the Kolmogorov bounds in \eqref{BE1}. 
	Here and thereafter, we denote $\Phi(z)=\gamma((-\infty,z])$ for $z\in\mathbb{R}.$
	
	\begin{proof}[Proof of the upper bound in \eqref{BE1}]
		To prove the upper bound in \eqref{BE1}, it suffices to consider $z>0$ since $W$ is symmetric.
		For $z>0$, we have
		\begin{equation}\label{pr00}
			|{\mathbb{P}}(W\le z)-\Phi(z)|\le  \max\{{\mathbb{P}}(W>z),1-\Phi(z)\}\le 0.5.
		\end{equation} If
		$N\le 100$, then the upper bound in \eqref{BE1} holds by \eqref{pr00}. Therefore we only need to consider $N>100$. For the case where $0<z<x_1$, there exists $1\le n\le m$ such that $z\in [x_{n+1},x_n)$. In this case, we have
		\[{\mathbb{P}}(W>z)=\dfrac{n}{N} ~\mbox{ and } ~ {\mathbb{P}}(W^*>z)=\dfrac{n-1}{N-1}+\varepsilon,\]
		where $\varepsilon\le 1/(N-1)$, and therefore
		\begin{equation*}
			\begin{split}
				\left|{\mathbb{P}}(W\le z)-{\mathbb{P}}(W^*\le z)\right|&=\left|\dfrac{n}{N}-\dfrac{n-1}{N-1}-\varepsilon\right|\\
				&\le \dfrac{n}{N}-\dfrac{n-1}{N-1}+\varepsilon\\
				&\le \dfrac{1}{N}+\dfrac{1}{N-1}\le \dfrac{2.02}{N},
			\end{split}
		\end{equation*}
		where we have applied $N>100$ in the last inequality. For the case where $z\ge x_1$, we have
		\[{\mathbb{P}}(W>z)={\mathbb{P}}(W^*>z)=0.\]
		Therefore,
		\begin{equation}
			\label{pr01}
			\begin{split}
				\sup_{z>0}\left|{\mathbb{P}}(W\le z)-{\mathbb{P}}(W^*\le z)\right|&\le \dfrac{2.02}{N}.
			\end{split}
		\end{equation}

		Now we bound the Kolmogorov distance between distribution of $W^*$ and $\gamma$. Let $z>0$,
		and let $f_{z}$ be
		the unique bounded solution of the Stein equation
		\[f^{'}(w)-w f(w)=1(w\le z)-\Phi(z),\]
		and 
		\[g_z(w)=(w f_z(w))^{'}.\]
		Then (see  Chen and Shao \cite[p. 248]{ChenShao01})
		\begin{equation}\label{ap01} g_z(w)=
			\begin{cases}
				\left(\sqrt{2\pi}(1+w^2)e^{w^2/2}(1-\Phi(w))-w\right)\Phi(z) &\mbox{ if  } w> z,\\
				\left(\sqrt{2\pi}(1+w^2)e^{w^2/2}\Phi(w)+w\right)(1-\Phi(z)) &\mbox{ if  } w\le z.\\
			\end{cases}
		\end{equation} 
		From Lemma 2.3 in Chen, Goldstein and Shao \cite{CGS}, we have $0<f_z(w)\le \sqrt{2\pi}/4$, $|f_{z}^{'}(w)|\le 1$, and
		therefore
		\begin{equation}\label{CS3a}
			|g_{z}(w)|\le |wf_{z}^{'}(w)|+|f_{z}(w)|\le |w|+\sqrt{2\pi}/4.
		\end{equation}
		Chen and Shao \cite[p. 249]{ChenShao01} proved that  $g_{z}\ge 0$, $g_{z}(w)\le
		2(1-\Phi(z))$ for $w\le 0$.
		
		We also have the following lemma whose proof is given in the Appendix.
		
		\begin{lem}\label{lem.gz}
			We have the following properties of $g_z$:
			\begin{equation}\label{CS3b}
				g_{z}(w) \text{ increases when }w\le z \text{ and decreases when }w>z,
			\end{equation}
			\begin{equation}\label{CS3c}
				g_{z}(w)\le \dfrac{3}{2(1-w)^3} \mbox{ for } w< 0,
			\end{equation}
			and
			\begin{equation}\label{CS3d}
				g_{z}(w)\le \dfrac{3}{(1+w)^3} \mbox{ for } w> z.
			\end{equation}
		\end{lem}
		Chen and Shao \cite{ChenShao01} proved that 
		\[ g_{z}(w)\le \dfrac{2}{1+w^3} \mbox{ for } w> z.\]
		For large $w$, this bound is of the same order as \eqref{CS3d} but with a better constant.
		We use \eqref{CS3d} because we need some technical
		estimates such as $(1+w)^3\ge x_{j-1}^3$ for $w\ge x_{j+1}$, $2\le j\le m-1$
		(see, for example, \eqref{obs3} and \eqref{pr13} below).
		
		The Kolmogorov distance 
		between the distribution of $W^*$ and $\gamma$ can be bounded as follows.
		\begin{equation}\label{pr02}
			\arraycolsep=1.2pt\def\arraystretch{2.6}
			\begin{array}{ll}
				&\left|{\mathbb{P}}(W^{*}\le z)-\Phi(z)\right|=\left|{\mathbb{E}}f_{z}^{'}(W^{*})-{\mathbb{E}}W^{*}f_{z}(W^{*})\right|\\
				&\quad =\left|\dfrac{N}{N-1}{\mathbb{E}}Wf_{z}(W)-{\mathbb{E}}W^{*}f_{z}(W^{*})\right|\\
				&\quad \le \left|{\mathbb{E}}Wf_{z}(W)-{\mathbb{E}}W^{*}f_{z}(W^{*})\right|+\dfrac{1}{N-1}{\mathbb{E}}|Wf_{z}(W)|.
			\end{array}
		\end{equation}
		Since $0<f_{z}(w)\le \sqrt{2\pi}/4$, the last term in \eqref{pr02} is bounded by
		\begin{equation}\label{pr03a}
			\dfrac{1}{N-1}{\mathbb{E}}|Wf_{z}(W)|\le \dfrac{\sqrt{2\pi {\mathbb{E}}W^2}}{4(N-1)}=\dfrac{\sqrt{2\pi}}{4 \sqrt{N(N-1)} }\le \dfrac{1}{N},
		\end{equation}
		where we have used $N>100$ in the last inequality. By the definition of $g_z$, we have
		\begin{equation}\label{pr03b}
			\begin{split}
				{\mathbb{E}}Wf_{z}(W)-{\mathbb{E}}W^{*}f_{z}(W^{*})= -{\mathbb{E}}\int_{0}^{W^*-W}g_{z}(W+t)dt.
			\end{split}
		\end{equation}
		For $n\le m$, from \eqref{lem31a} of Lemma \ref{3} below, we see that $x_n$ is bounded by an absolute constant as long as $\log(m/n)$
		is bounded. 
		To bound the right hand side of \eqref{pr03b}, we separate
		$W$ around $m/e^3$ (we can replace $e^3$ by any constant $c>1$, the number $e^3$ is only for convenience 
		when we compute the explicit constant in \eqref{BE1}) as follows.
		\begin{equation}\label{pr03b1}
			\begin{split}
				R_1= {\mathbb{E}}\left|\int_{0}^{W^*-W}g_{z}(W+t){\mathbf{1}}\left(W\le  -x_{\lfloor m/e^3\rfloor}\right) dt\right|,
			\end{split}
		\end{equation}
		\begin{equation}\label{pr03b2}
			\begin{split}
				R_2= {\mathbb{E}}\left|\int_{0}^{W^*-W}g_{z}(W+t){\mathbf{1}}\left(|W|< x_{\lfloor m/e^3\rfloor}\right) dt\right|,
			\end{split}
		\end{equation}
		and
		\begin{equation}\label{pr03b3}
			\begin{split}
				R_3= {\mathbb{E}}\left|\int_{0}^{W^*-W}g_{z}(W+t){\mathbf{1}}\left(W\ge  x_{\lfloor m/e^3\rfloor}\right) dt\right|.
			\end{split}
		\end{equation}

	We need the following lemma whose proof will be presented 
		in the Appendix. 
		\begin{lem}\label{3} Let $N>100$. The following statements hold.
			\begin{equation}\label{lem31d}
				0\le x_m \le \frac{1}{m}.
			\end{equation}
			For $1\le n\le m-1$, we have
			\begin{equation}\label{lem31a}
				x_n \le \sqrt{2(1+\log(m/n))},
			\end{equation}
			and
			\begin{equation}\label{lem31b}
				\left(\frac{n(n+1)}{2}\right)^{1/2}\le S_n \le \dfrac{3n}{2}\sqrt{2(1+\log(m/n))}.
			\end{equation}
			For $1\le n\le \lfloor m/e^3\rfloor$, we have
			\begin{equation}\label{lem31c}
				x_n\ge \dfrac{1}{3}\sqrt{2(1+\log(m/n))}~\text{and}~S_n\ge \dfrac{n}{3}\sqrt{2(1+\log(m/n))}.
			\end{equation}
			For $1\le i<j\le \lfloor m/e^3\rfloor$, we have
			\begin{equation}\label{lem31e}
				x_{i}^2-x_{j}^2\ge  \dfrac{4\log(j/i)}{9}.
			\end{equation}
		\end{lem}

		\begin{rem}\label{2}
			\rm{ 
				By using the first half of \eqref{lem31c}, we have $x_{\lfloor m/e^3\rfloor}\ge \sqrt{8}/3$.
				We note also that $m>50$ and $S_j\ge jx_j$, $1\le j\le m$.
				These simple inequalities will be used in many places later without further mention.
				For $x_1$, McKeague and Levin \cite{ML16} proved the following lower bound
				\begin{equation}\label{ML02}
					S_1=x_1\ge \sqrt{\log(m)},
				\end{equation}
				which is of the same order as ours but with a better constant, but their method seems
				not to work with $x_n$ for $n>1$.
			}
		\end{rem}

		We will now bound $R_1$, $R_2$ and $R_3$. 
		From \eqref{MIW-recursion} and the first half of \eqref{lem31b},
		we have $x_1-x_2=1/S_1\le 1,$
		and
		\[x_{n-1}-x_{n+1}=\dfrac{1}{S_{n-1}}+\dfrac{1}{S_n}\le \left(\dfrac{2}{(n-1)n}\right)^{1/2}+\left(\dfrac{2}{n(n+1)}\right)^{1/2}\le 1\]
		for $3\le n\le m$, or equivalently, 
		\begin{equation}\label{obs3}
			x_1\le x_2+1\text{ and } x_{n-1}\le x_{n+1}+1\text{ for }3\le n\le m.
		\end{equation}
		Let $t$ be a real number lying between $0$ and $W^*-W$.
		From the definitions of $W$ and $W^*$, we observe the following facts.
		\begin{obs}\label{obs1}
			\begin{description}
				\item{(i)} If $W=x_1$, then 
				\[|W^*-W|\le x_1-x_2\ \text{ and }\ x_2\le  W+t\le x_1;\]
				\item{(ii)} If $W=x_n$, $2\le n\le m$, then 
				\[|W^*-W|\le x_{n-1}-x_n\ \text{ and }\ x_{n+1}\le W+t\le x_{n-1};\]
			\end{description}
			and, by symmetry, we have
			\begin{description}
				\item{(iii)} If $W=x_N=-x_1$, then 
				\[|W^*-W|\le x_1-x_2\ \text{ and }\ -x_1\le W+t\le -x_2;\]
				\item{(iv)} If $W=x_{N+1-n}=-x_n$, $2\le n\le m$, then 
				\[|W^*-W|\le x_{n-1}-x_n\ \text{ and }\ -x_{n-1}\le W+t\le -x_{n+1}.\]
			\end{description}
		\end{obs}
		Keeping Facts \ref{obs1} and the properties of $g_z$ in mind, we have
		\begin{equation}\label{pr13}
			\begin{split}
				R_{1}&\le \dfrac{1}{N}\left((x_1-x_2)g_{z}(-x_2)+\sum_{j=2}^{\lfloor m/e^3\rfloor}(x_{j-1}-x_{j})g_{z}(-x_{j+1})\right)\\
				&\le \dfrac{3}{2N}\left(\dfrac{x_1-x_2}{(1+x_2)^3}+\sum_{j=2}^{\lfloor m/e^3\rfloor}\dfrac{x_{j-1}-x_{j}}{(1+x_{j+1})^3}\right)\\
				&= \dfrac{3}{2N}\left(\dfrac{1}{x_{1}(1+x_2)^3}+\sum_{j=2}^{\lfloor m/e^3\rfloor}\dfrac{1}{(1+x_{j+1})^3S_{j-1}}\right)\\
				&\le\dfrac{3}{2N}\left(\dfrac{1}{x_{1}^4}+\sum_{j=2}^{\lfloor m/e^3\rfloor}\dfrac{1}{(j-1)x_{j-1}^4}\right)\\
				&\le\dfrac{3}{2N}\left(\dfrac{2}{x_{1}^4}+\sum_{j=2}^{\lfloor m/e^3\rfloor-1}\dfrac{81}{4j(1+\log(m/j))^2}\right)\\
				&\le \dfrac{3}{2N}\left(\dfrac{2}{\log^2(m)}+\int_{1}^{m/e^3}\dfrac{81dx}{4x(1+\log(m/x))^2}\right)\\
				&= \dfrac{3}{2N}\left(\dfrac{2}{\log^2(m)}+\dfrac{81}{16}-\dfrac{81}{4(1+\log(m))}\right)\le \dfrac{243}{32N},
			\end{split}
		\end{equation}
		where we have applied \eqref{CS3b} in the first inequality, and \eqref{CS3c} in the second inequality, \eqref{obs3} and simple bounds $S_j\ge jx_j$ ($1\le j\le m$) in the 
		third inequality,
		the first half of \eqref{lem31c} and \eqref{ML02} in the fourth inequality. The final bound in \eqref{pr13}
		follows from $m>50$ and the elementary inequality $2/\log^2(m)\le 81/(4(1+\log(m)))$.
		For $R_2$, we have
		\begin{equation}\label{pr14}
			\begin{split}
				R_{2}&\le \dfrac{2}{N}\sum_{0\le x_j<\lfloor m/e^3\rfloor}(x_{j-1}-x_{j})(x_{j-1}+\sqrt{2\pi}/4)\\
				&=\dfrac{2}{N}\sum_{\lfloor m/e^3\rfloor< j\le m}\left(\dfrac{x_{j-1}}{S_{j-1}}+\dfrac{\sqrt{2\pi}}{4S_{j-1}}\right)\\
				&\le\dfrac{2}{N}\sum_{\lfloor m/e^3\rfloor< j\le m}\left(\dfrac{1}{j-1}+\dfrac{\sqrt{\pi}}{2(j-1)}\right)\\
				&\le \dfrac{2+\sqrt{\pi}}{N}\left(\dfrac{1}{\lfloor m/e^3\rfloor}+\dfrac{1}{\lfloor m/e^3\rfloor+1}+\int_{m/e^3}^m\dfrac{dx}{x}\right)\\
				&\le \dfrac{2+\sqrt{\pi}}{N}\left(\dfrac{1}{2}+\dfrac{1}{3}+3\right)\le \dfrac{15}{N},
			\end{split}
		\end{equation}
		where we have applied \eqref{CS3a} and Facts \ref{obs1} in the first inequality, 
		the first half of \eqref{lem31b} and simple bounds $x_j/S_j\le 1/j$ ($1\le j\le m$) in the second inequality,
		and used $m>50$ in the fourth inequality.
		
		To bound $R_3$, we set $\ell=\min\{j: x_j\le z\}$ and consider the following two cases.
		
		{\it Case 1}: $x_{\ell}<x_{\lfloor m/e^3\rfloor+1}$. In this case, we have $z< x_{\lfloor m/e^3\rfloor+1}$. If $W=x_j\ge x_{\lfloor m/e^3\rfloor}$, then, by Facts \ref{obs1}, 
		\begin{equation}\label{pr15}
			W+t>x_{j+1} >z\ \text{ for all }\ t \ \text{ lying between $0$ and } W^*-W.
		\end{equation}
		Therefore,
		\begin{equation}\label{boundR3case1}
			\begin{split}
				R_{3}& = {\mathbb{E}}\left|\int_{0}^{W^*-W}g_{z}(W+t){\mathbf{1}}\left(W\ge x_{\lfloor m/e^3\rfloor}\right) dt\right|\\
				&\le \dfrac{1}{N}\left((x_1-x_2)g_{z}(x_2)+\sum_{j=2}^{\lfloor m/e^3\rfloor}(x_{j-1}-x_{j})g_{z}(x_{j+1})\right)\\
				&\le \dfrac{3}{N}\left(\dfrac{x_1-x_2}{(1+x_2)^3}+\sum_{j=2}^{\lfloor m/e^3\rfloor}\dfrac{x_{j-1}-x_{j}}{(1+x_{j+1})^3}\right)\le \dfrac{243}{16N},
			\end{split}
		\end{equation}
		where we have applied \eqref{CS3b} and \eqref{pr15} in the first inequality,
		and \eqref{CS3d} and \eqref{pr15} in the second inequality. The last inequality in \eqref{boundR3case1} follows by using 
		the same calculations as in \eqref{pr13}.
		
		{\it Case 2}: $x_{\ell}\ge x_{\lfloor m/e^3\rfloor+1}$.  
		In this case, we have $z\ge x_{\lfloor m/e^3\rfloor+1}$, and $R_3$ is bounded by $R_{31}+R_{32}+R_{33}$,
		where
		\begin{equation}\label{pr16}
			\begin{split}
				R_{31}& = {\mathbb{E}}\left|\int_{0}^{W^*-W}g_{z}(W+t){\mathbf{1}}\left(x_{\lfloor m/e^3\rfloor}\le W<x_{\ell}\right) dt\right|,
			\end{split}
		\end{equation}
		\begin{equation}\label{pr16b}
			\begin{split}
				R_{32}& = {\mathbb{E}}\left|\int_{0}^{W^*-W}g_{z}(W+t){\mathbf{1}}\left(x_{\ell}\le W\le x_{\ell-1}\right) dt\right|,
			\end{split}
		\end{equation}
		and
		\begin{equation}\label{pr16c}
			\begin{split}
				R_{33}& = {\mathbb{E}}\left|\int_{0}^{W^*-W}g_{z}(W+t){\mathbf{1}}\left(W>x_{\ell-1}\right) dt\right|.
			\end{split}
		\end{equation}
		Here and thereafter, we denote $x_0=x_1$ and $S_0=S_1$.
		We of course only need to bound $R_{31}$ when $x_{\lfloor m/e^3\rfloor}< x_{\ell}$.  
		Since $x_{\ell}\le z$, 
		we conclude that if $x_{\lfloor m/e^3\rfloor}\le W=x_{j+1}<x_{\ell}$, then by Facts \ref{obs1}, we have
		\begin{equation}\label{pr16d}
			|W^*-W|\le x_j-x_{j+1} \text{ and }	0<W+t<x_{j}\le x_{\ell} \le z
		\end{equation}
		for all $t$ lying between $0$ and $W^*-W$.
		Applying \eqref{ap01}, \eqref{pr16d} and Facts \ref{obs1}, we have
		\begin{equation}\label{pr18}
			\begin{split}
				R_{31}&\le \dfrac{1-\Phi(x_{\ell})}{N}\left(\sum_{j={\ell}}^{\lfloor m/e^3\rfloor-1}(x_{j}-x_{j+1}) \left(\sqrt{2\pi}(1+x_{j}^2)e^{x_{j}^2/2}+x_j\right)\right).
			\end{split}
		\end{equation}
		Since $1-\Phi(x_{\ell})\le e^{-x_{\ell}^2/2}/\left(x_{\ell}\sqrt{2\pi}\right),$ 
		it follows from \eqref{pr18} that
		\begin{equation}\label{pr18a}
			\begin{split}
				R_{31}&\le \dfrac{1}{N}\sum_{j={\ell}}^{\lfloor m/e^3\rfloor-1}\left(\dfrac{1+x_{j}^2}{S_jx_{\ell}}e^{(x_{j}^2-x_{\ell}^2)/2}+\dfrac{x_j e^{-x_{\ell}^2/2}}{\sqrt{2\pi}S_jx_{\ell}}\right).
			\end{split}
		\end{equation}
		For all $\ell\le j\le \lfloor m/e^3\rfloor$, by applying \eqref{lem31c}, we have
		\begin{equation}\label{pr18b}
			S_jx_{\ell}\ge S_jx_j\ge 8j/9.
		\end{equation}
		By using \eqref{lem31e}, \eqref{pr18b}, and simple inequalities $x_{\ell}\ge x_{\lfloor m/e^3\rfloor}\ge \sqrt{8}/3$ and $x_j/S_j\le 1/j$ ($1\le j\le m$), we have
		\begin{equation}\label{pr18c}
			\begin{split}
				&\dfrac{1+x_{j}^2}{S_jx_{\ell}}e^{(x_{j}^2-x_{\ell}^2)/2}+\dfrac{x_j e^{-x_{\ell}^2/2}}{\sqrt{2\pi}S_jx_{\ell}}\\
				&\le \left(\dfrac{1}{S_jx_j}+\dfrac{x_j}{S_j}\right)\left(\dfrac{\ell}{j}\right)^{2/9}+\dfrac{3x_j}{4\sqrt{\pi}S_j}\exp\left(-\dfrac{1}{2}\left(\dfrac{\sqrt{2\log(m/{\ell})}}{3}\right)^2\right)\\
				&\le \left(\dfrac{9}{8j}+\dfrac{1}{j}\right)\left(\dfrac{\ell}{j}\right)^{2/9}+\dfrac{3}{4\sqrt{\pi}j}\left(\dfrac{\ell}{m}\right)^{1/9}
			\end{split}
		\end{equation}
		for all $\ell\le j\le \lfloor m/e^3\rfloor$. It follows from \eqref{pr18a} and \eqref{pr18c} that
		\begin{equation}\label{pr19}
			\begin{split}
				&R_{31}\le \dfrac{1}{N}\left(\dfrac{17}{8}\sum_{j={\ell}}^{\lfloor m/e^3\rfloor-1}\dfrac{{\ell}^{2/9}}{j^{11/9}}+\dfrac{3e^{-3}}{4\sqrt{\pi}}\sum_{j={\ell}}^{\lfloor m/e^3\rfloor-1}\dfrac{{\ell}^{1/9}}{j^{10/9}}\right)\\
				&\le \dfrac{1}{N}\left(\dfrac{17}{8}\left(\dfrac{1}{\ell}+ \int_{\ell}^{m/e^3}\dfrac{{\ell}^{2/9}dx}{x^{11/9}}\right)+\dfrac{3e^{-3}}{4\sqrt{\pi}}\left(\dfrac{1}{\ell}+\int_{\ell}^{m/e^3}\dfrac{{\ell}^{1/9}dx}{x^{10/9}}\right)\right)\\
				&\le\dfrac{1}{N}\left(\dfrac{17}{8}\left(1+\dfrac{9}{2}\right)+\dfrac{30e^{-3}}{4\sqrt{\pi}}\right)\le\dfrac{12}{N}.
			\end{split}
		\end{equation}
		
		Now, we bound $R_{32}$ and $R_{33}$. 
		If $x_{\ell}\ge x_3$, then 
		\begin{equation}\label{pr17}
			\begin{split}
				R_{32}&+R_{33} \le \sum_{i=1}^3{\mathbb{E}}\left|\int_{0}^{W^*-W}g_z\left(W+t\right){\mathbf{1}}\left(W= x_{i}\right) dt\right|\\
				&\le \dfrac{1}{N}\left(2(x_{1}-x_{2})\left(x_{1}+\dfrac{\sqrt{2\pi}}{4}\right)+(x_{2}-x_{3})\left(x_{2}+\dfrac{\sqrt{2\pi}}{4}\right)\right)\\
				&= \dfrac{1}{N}\left(\dfrac{2(x_{1}+\sqrt{2\pi}/4)}{S_{1}}+\dfrac{x_{2}+\sqrt{2\pi}/4}{S_{2}}\right)\le\dfrac{6}{N},
			\end{split}
		\end{equation}
		where we have applied \eqref{CS3a} and Facts \ref{obs1} in the second inequality,
		and the first half of \eqref{lem31b} and a simple bound $x_2/S_2\le 1/2$ in the last inequality. 
		If $x_{\ell}\le x_4$, then, similarly to \eqref{pr17}, we have
		\begin{equation}\label{pr20}
			\begin{split}
				R_{32}& = {\mathbb{E}}\left|\int_{0}^{W^*-W}g_z\left(W+t\right){\mathbf{1}}\left(x_{\ell}\le W\le x_{\ell-1}\right) dt\right|\\
				&\le \dfrac{1}{N}\left(\dfrac{x_{\ell-1}+\sqrt{2\pi}/4}{S_{\ell-1}}+\dfrac{x_{\ell-2}+\sqrt{2\pi}/4}{S_{\ell-2}}\right)\\
				&\le \dfrac{1}{N}\left(\dfrac{x_3}{S_3}+\dfrac{\sqrt{2\pi}/4}{S_3}+\dfrac{x_2}{S_2}+\dfrac{\sqrt{2\pi}/4}{S_2}\right)\\
				&\le \dfrac{1}{N}\left(\dfrac{1}{3}+\dfrac{\sqrt{2\pi}/4}{S_3}+\dfrac{1}{2}+\dfrac{\sqrt{2\pi}/4}{S_2}\right)\le \dfrac{2}{N},
			\end{split}
		\end{equation}
		and
		\begin{equation}\label{pr21}
			\begin{split}
				R_{33}&= {\mathbb{E}}\left|\int_{0}^{W^*-W}g_z(W+t){\mathbf{1}}\left(W>x_{\ell-1}\right) dt\right|\\
				&\le \dfrac{3}{N}\left(\dfrac{x_{1}-x_{2}}{(1+x_{2})^3}+\sum_{j=2}^{l-2}\dfrac{x_{j-1}-x_{j}}{(1+x_{j+1})^3}\right)\\
				&\le \dfrac{243}{16N},
			\end{split}
		\end{equation}
		where, in \eqref{pr21}, we have applied Facts \ref{obs1} and \eqref{CS3d}
		in the first inequality, and used the same calculations as in \eqref{pr13} in the second inequality.
		Combining \eqref{pr19}--\eqref{pr21}, we have
		\begin{equation}\label{pr23}
			R_3\le \dfrac{467}{16N}.
		\end{equation}
		It follows from \eqref{pr13}, \eqref{pr14} and \eqref{pr23} that
		\begin{equation}\label{pr24}
			\left|{\mathbb{E}}Wf_{z}(W)-{\mathbb{E}}W^{*}f_{z}(W^{*})\right|\le R_1+R_2+R_3\le \dfrac{1657}{32N}.
		\end{equation}
		Combining \eqref{pr02}, \eqref{pr03a}, and \eqref{pr24}, we have
		\begin{equation}\label{pr25}
			\begin{split}
				|{\mathbb{P}}(W^*\le z)-\Phi(z)|&\le \dfrac{1}{N}+\dfrac{1657}{32N}= \dfrac{1689}{32N}.
			\end{split}
		\end{equation}
		The upper bound of \eqref{BE1} follows from \eqref{pr01} and \eqref{pr25}.
	\end{proof}
	
	\begin{proof}[Proof of the lower bound in \eqref{BE1}]
		The proof of the lower bound in \eqref{BE1} follows from
		the fact that the Kolmogorov distance between $\gamma$ and the
		probability distribution of a discrete random variable is 
		always greater than half of the minimum of the jumps. In our case, the jumps are all equal to $1/N$. Therefore,
		$$d_{\mathrm{K}}({\mathbb{P}}_N,\gamma)\ge \dfrac{1}{2N}.$$ 
	\end{proof}

	\section{Bounds on the Wasserstein distance}\label{sec:Wasserstein}
	
	\begin{proof}[Proof of \eqref{BE2}]
		
		By Theorem 1.1 of Goldstein  \cite{Goldstein07}, we have
		\begin{equation}
			\label{lp00}
			\begin{split} 
				d_{\mathrm{W}}({\mathbb{P}}_N,\gamma)&\le 2{\mathbb{E}}|W-W^*|\\
				&\le\dfrac{2}{N-1}\sum_{n=1}^{N-1}(x_n-x_{n+1})\\
				&=\dfrac{4x_1}{N-1}\le \dfrac{4\sqrt{2(1+\log(m))}}{N-1} \\
				&\le \dfrac{8\sqrt{\log N}}{N-1}\le \dfrac{16\sqrt{\log N}}{N},
			\end{split}
		\end{equation}
		where we have applied \eqref{lem31a} in the third inequality, and definition of $m$ 
		in the fourth inequality. The proof of the upper bound of \eqref{BE2} is completed.
		
		For the lower bound of \eqref{BE2}, we can assume $N>100$ since it is trivial when $N\le 100$ and $C\ge 3$. We note that $d_{\mathrm{W}}({\mathbb{P}}_N,\gamma) \ge \left|{\mathbb{E}}h(W)-{\mathbb{E}}h(Z)\right|$ for any $1$-Lipschitz function $h$. We will use the ``sawtooth'' piecewise linear function considered in McKeague and Levin \cite{ML16}, which is defined as follows:
		\begin{equation*}
			h(w)=
			\begin{cases}
				0& \mbox{ if }w>x_1 \mbox{ or } w<x_N,\\
				w-x_{n+1}& \mbox{ if } x_{n+1}\le w <m_n,\ 1\le n< N,\\
				x_n-w& \mbox{ if } m_n\le w<x_n,\ 1\le n< N,
			\end{cases}
		\end{equation*}
where $m_n = (x_n +x_{n+1})/2, 1\le n< N$. 
		Clearly $h$ is $1$-Lipschitz and ${\mathbb{E}}h(W) = 0$. Simple calculation gives 
		\begin{equation}\label{lp01}
			{\mathbb{E}}h(W^*) = \dfrac{x_1}{2(N-1)}\ge \dfrac{\sqrt{\log (N/2)}}{2(N-1)}\ge \dfrac{\sqrt{\log (N/2)}}{2N},
		\end{equation}
		where we have applied \eqref{ML02} in the first inequality.
		
		For this ``sawtooth'' function $h$,
		let $f_{h}$ be
		the unique bounded solution of the Stein equation
		\[f^{'}(w)-w f(w)=h(w)-{\mathbb{E}}h(Z),\]
		and let
		\[g_h(w)=(w f_{h}(w))^{'}.\]
		By Lemma 2.4 in \cite{CGS}, we have $\|f_h\|\le 2$ and $\|f_{h}^{'}\|\le \sqrt{2/\pi}$.
		Therefore, 
		\begin{equation}\label{g01}
			|g_h(w)|\le |f_{h}(w)|+|wf^{'}(w)|\le 2(1+|w|).
		\end{equation}
		It thus follows from \eqref{lp01} that
		\begin{equation}\label{lp02}
			\begin{split} 
				&\left|{\mathbb{E}}h(W)-{\mathbb{E}}h(Z)\right| = \left|{\mathbb{E}}f_{h}^{'}(W)-{\mathbb{E}}Wf_{h}(W)\right|\\
				&\qquad=\left|{\mathbb{E}}f_{h}^{'}(W)-\left(1-\dfrac{1}{N}\right){\mathbb{E}}f_{h}^{'}(W^{*})\right|\\
				&\qquad=\left|{\mathbb{E}}Wf_{h}(W)-{\mathbb{E}}W^*f_{h}(W^{*})-{\mathbb{E}}h(W^*)+\dfrac{1}{N}{\mathbb{E}}f_{h}^{'}(W^*)\right|\\
				&\qquad\ge \left|{\mathbb{E}}h(W^*)\right|-\left|{\mathbb{E}}Wf_{h}(W)-{\mathbb{E}}W^*f_{h}(W^{*})\right|-\left|\dfrac{1}{N}{\mathbb{E}}f_{h}^{'}(W^*)\right|\\
				&\qquad\ge \dfrac{\sqrt{\log (N/2)}}{2N}-\left|{\mathbb{E}}Wf_{h}(W)-{\mathbb{E}}W^*f_{h}(W^{*})\right|-\dfrac{\sqrt{2/\pi}}{N}.
			\end{split}
		\end{equation}
		The lower bound  of \eqref{BE2} will follow if we can show that
		\[\left|{\mathbb{E}}Wf_{h}(W)-{\mathbb{E}}W^*f_{h}(W^{*})\right|\le \dfrac{C}{N}.\]
		Letting $T=W^*-W$, we have
		\begin{equation}\label{lp10}
			{\mathbb{E}}|Wf_{h}(W)-W^*f_{h}(W^*)|={\mathbb{E}}\left(|Tg_h(W+\xi)|\right),
		\end{equation}
		where $\xi$ is a random variable lying between $0$ and $T$. We will obtain an upper bound on ${\mathbb{E}}\left(|Tg_h(W+\xi)|\right)$ using the the same truncation for $W$ as as in the proof for the Kolmogorov bound. So we let
		
		\begin{equation}\label{r1}
			R_1 = {\mathbb{E}}\left(|Tg_h(W+\xi)|{\mathbf{1}}(|W|\le x_{\lfloor m/e^3\rfloor})\right),
		\end{equation}
		and
		\begin{equation}\label{r2}
			R_2 = {\mathbb{E}}\left(|Tg_h(W+\xi)|{\mathbf{1}}(|W|> x_{\lfloor m/e^3\rfloor})\right).
		\end{equation}
		For $R_1$, we have
		\begin{equation}\label{lp12add01}
			\begin{split} 
				R_1 &\le 2{\mathbb{E}}\left(|T|(1+|W+\xi|){\mathbf{1}}(|W|\le x_{\lfloor m/e^3\rfloor})\right)\\
				&\le \dfrac{C}{N}\sum_{i=\lfloor m/e^3\rfloor}^{m}(x_{i-1}-x_i)(1+x_{i-1})\\
				&= \dfrac{C}{N}\sum_{i=\lfloor m/e^3\rfloor}^{m} \left(\dfrac{1}{S_{i-1}}+\dfrac{x_{i-1}}{S_{i-1}}\right)\\
				&\le \dfrac{C}{N}\sum_{i=\lfloor m/e^3\rfloor}^{m} \left(\dfrac{1}{S_{i-1}}+\dfrac{1}{i-1}\right),
			\end{split}
		\end{equation}
		where we have applied \eqref{g01} in the first inequality, Facts \ref{obs1} in the second inequality, and simple bounds $x_j/S_j\le 1/j$ ($1\le j\le m$) in the third inequality.
		If $\lfloor m/e^3\rfloor\le i\le m$, then by using the second half of \eqref{lem31c}, we have $S_{i-1}\ge S_{\lfloor m/e^3\rfloor-1}\ge m/C$. Therefore, \eqref{lp12add01} implies
		\begin{equation}\label{lp12}
			\begin{split} 
				R_1 &\le \dfrac{C}{N}\sum_{i=\lfloor m/e^3\rfloor}^{m} \left(\dfrac{1}{m}+\dfrac{1}{i-1}\right)\le \dfrac{C}{N}.
			\end{split}
		\end{equation}
		To bound $R_2$, we need the following lemma.
		\begin{lem}\label{lem41}
			If either $x_{\lfloor m/e^3\rfloor}\le x_{n+1}< w \le x_n$ or $-x_n\le w< -x_{n+1}\le -x_{\lfloor m/e^3\rfloor}$, then
			\begin{equation}\label{lem41_01}
				\begin{split} 
					&|g_h(w)|\le  C\left(\dfrac{1}{\log(m/n)}+\dfrac{1}{n^{2/9}}\right).
				\end{split}
			\end{equation}
		\end{lem}
		
		We postpone the proof of \lemref{lem41} to the Appendix.
		Now, we bound $R_2$ as follows.
		\begin{equation}\label{lp13}
			\begin{split}
				R_2 &\le \dfrac{C}{N}\sum_{n=2}^{\lfloor m/e^3\rfloor-1}(x_{n-1}-x_{n})\left(\dfrac{1}{\log(m/n)}+\dfrac{1}{n^{2/9}}\right)\\
				&= \dfrac{C}{N}\sum_{n=2}^{\lfloor m/e^3\rfloor-1}\left(\dfrac{1}{\log(m/n)S_{n-1}}+\dfrac{1}{n^{2/9}S_{n-1}}\right)\\
				&\le \dfrac{C}{N}\sum_{n=1}^{\lfloor m/e^3\rfloor-2}\left(\dfrac{1}{n\log^{3/2}(m/n)}+\dfrac{1}{n^{11/9}\log^{1/2}(m/n)}\right)\\
				&\le \dfrac{C}{N},
			\end{split}
		\end{equation}
		where we have applied Facts \ref{obs1} and Lemma \ref{lem41} in the first inequality, and the second half of \eqref{lem31c}
		in the second inequality. 
		It follows from \eqref{lp12} and \eqref{lp13} that 
		\begin{equation}\label{lp15}
			\begin{split} 
				{\mathbb{E}}|Wf_{h}(W)-W^*f_{h}(W^*)|& \le \dfrac{C}{N}.
			\end{split}
		\end{equation}
		Combining  \eqref{lp02} and \eqref{lp15}, we have
		\begin{equation*}
			\label{lp16}
			d_{\mathrm{W}}({\mathbb{P}}_N,\gamma)\ge \dfrac{\sqrt{\log (N/2)}}{2N}-\dfrac{C}{N},
		\end{equation*}
		thereby proving the lower bound of \eqref{BE2}.
	\end{proof}

	\appendix
	\section{}\label{sec:Appendix}
	In this appendix, we will prove those lemmas in the previous sections, whose proofs have been deferred to the appendix.
	
	\begin{proof}[Proof of \lemref{lem.gz}]
		We first prove that
		\begin{equation}\label{ap02}
			\sqrt{2\pi} we^{w^2/2}(1-\Phi(w))\le \dfrac{w^2+2}{w^2+3}
		\end{equation}
		for $w> 0$.
		This inequality is equivalent to
		\begin{equation}\label{ap04}
			\int_{w}^{\infty}e^{-x^2/2}dx\le \dfrac{w^2+2}{w^3+3w}e^{-w^2/2}.
		\end{equation}
		Let $$h(w)=\dfrac{w^2+2}{w^3+3w}e^{-w^2/2}.$$
		Then
		$$h(w)=\int_{w}^\infty (-h'(x))dx.$$
		For $x>0$, we have 
		$$-h'(x)=\left(1+\dfrac{6}{(x^3+3x)^2}\right)e^{-x^2/2}\ge e^{-x^2/2}$$
		thereby proving \eqref{ap04}. The proof of \eqref{ap02} is completed.
		
		From \eqref{ap01} and \eqref{ap02}, we have for $w>z$, 
		$$g_{z}^{'}(w)=\left(\sqrt{2\pi}(w^2+3)we^{w^2/2}(1-\Phi(w))-2-w^2\right)\Phi(z)\le 0,$$
		which proves $g_z(w)$ is decreasing for $w>z$. Similarly, 
		$g_z(w)$ is increasing for $w\le z$ since in this case,
		$$g_{z}^{'}(w)=\left(2+w^2-\sqrt{2\pi}(w^2+3)(-w)e^{w^2/2}(1-\Phi(-w))\right)\left(1-\Phi(z)\right)\ge 0.$$
		The proof of \eqref{CS3b} is completed.
		
		Now we prove that for $w> 0$,
		\begin{equation}\label{ap06}
			0<\sqrt{2\pi}(1+w^2)e^{w^2/2}(1-\Phi(w))-w \le \dfrac{3}{(1+w)^3}.
		\end{equation}
		The first inequality in \eqref{ap06} is proved by Chen and Shao \cite{ChenShao01}.
		The second inequality in \eqref{ap06} is equivalent to
		\begin{equation}\label{ap08}
			\int_{w}^{\infty}e^{-x^2/2}dx\le \left(\dfrac{w}{1+w^2}+\dfrac{3}{(1+w^2)(1+w)^3}\right) e^{-w^2/2}.
		\end{equation}
		Let $$k(w)=\left(\dfrac{w}{1+w^2}+\dfrac{3}{(1+w^2)(1+w)^3}\right) e^{-w^2/2},\ w>0.$$
		Then
		$$k(w)=\int_{w}^\infty (-k'(x))dx.$$
		For $x>0$, we have 
		$$-k'(x)=\left(1+\dfrac{x^4-5x^3+6x^2+x+7}{(1+x^2)^2(1+x)^4}\right)e^{-x^2/2}\ge e^{-x^2/2}$$
		thereby proving \eqref{ap08}. Combining \eqref{ap01} and \eqref{ap06}, and noting that 
		$1-\Phi(z)\le 1/2,$ we have
		$$0\le g_z(w)\le \dfrac{3}{(1+w)^3} \text{ for }w>z,$$
		and 
		$$0\le g_z(w)\le \dfrac{3}{2(1-w)^3} \text{ for }w<0.$$
		This ends the proof of \eqref{CS3c} and \eqref{CS3d}.
	\end{proof}

	\begin{proof}[Proof of \lemref{3}]
		From (P3) and (P4) of \lemref{ML1}, it is easy to see that
		$x_m=0$ if  $N$ is odd, and $2x_m=x_m-x_{m+1}>0$ if $N$ is even.
		So we only need to prove the second inequality in \eqref{lem31d} for the case where $N$ is even.
		By \eqref{MIW-recursion}, we have
		$$x_j=\sum_{i=j}^{m-1}(x_i-x_{i+1})+x_m\ge \sum_{i=j}^{m-1}\dfrac{1}{S_i},\ 1\le j\le m-1.$$
		Therefore, $S_1\ge 1$ and
		\begin{equation}\label{21}
			\begin{split} S_j & \ge\dfrac{1}{S_1}+\dfrac{2}{S_2}+\dots+\dfrac{j-1}{S_{j-1}}+j\sum_{i=j}^{m-1}\dfrac{1}{S_i}\\
				& \ge \dfrac{1}{S_j}+\dfrac{2}{S_j}+\dots+\dfrac{j}{S_j} =\dfrac{j(j+1)}{2S_j},
			\end{split}
		\end{equation}
		for $2\le j\le m-1$.
		This implies that 
		\begin{equation}\label{22}
			S_j \ge \left(\frac{j(j+1)}{2}\right)^{1/2},~ 1\le j\le m-1.
		\end{equation}
		Since $x_m\ge 0$, $$S_m\ge S_{m-1}\ge \left(\dfrac{m(m-1)}{2}\right)^{1/2}.$$ 
		If $N$ is even, then $x_m = - x_{m+1}$. So
		$$2x_m=x_m-x_{m+1}=\dfrac{1}{S_m}\le \left(\dfrac{2}{m(m-1)}\right)^{1/2}.$$
		Therefore, 
		$$0\le x_m\le \left(\dfrac{1}{2m(m-1)}\right)^{1/2}\le \dfrac{1}{m}.$$
		This proves \eqref{lem31d}.
		
		Since $\{x_1,\dots,x_N\}$ is decreasing, $jx_j\le S_j\le jx_1$. For $1\le n\le m-1$, by using the telescoping sum and \eqref{lem31d}, we have
		\begin{equation}\label{lem35}
			\begin{split}
				x_{n}^2 & =\sum_{j=n}^{m-1}(x_{j}^2-x_{j+1}^2)+x_{m}^2 \le 2\sum_{j=n}^{m-1}x_j(x_{j} -x_{j+1})+x_{m}^2\\
				& =2\sum_{j=n}^{m-1}\dfrac{x_j}{S_j}+x_{m}^2 \le 2\sum_{j=n}^{m}\dfrac{1}{j}\\
				& \le 2(1/n+\log(m/n))\\
				& \le 2(1+\log(m/n)).
			\end{split}
		\end{equation} 
		This proves \eqref{lem31a}. For $2\le n\le m-1$,
		\begin{equation}\label{cal01}
			\begin{split} S_n &=x_1+\dots+x_n\\
				&\le \sqrt{2}\left(\sqrt{1+\log(m)}+\sqrt{1+\log(m/2)}+\dots+\sqrt{1+\log(m/n)}\right)\\
				&\le \sqrt{2}\left(\sqrt{1+\log(m)}+\int_{1}^{n}\sqrt{1+\log(m/x)}dx\right).
			\end{split}
		\end{equation}
		Set 
		\[I_n=\int_{1}^{n}\sqrt{1+\log(m/x)}dx.\] By using integration by parts, we have
		\begin{equation}\label{cal02}
			\begin{split} I_n & =n\sqrt{1+\log(m/n)}-\sqrt{1+\log(m)}+\dfrac{1}{2}\int_{1}^{n}\dfrac{dx}{\sqrt{1+\log(m/x)}}\\
				& \le n\sqrt{1+\log(m/n)}-\sqrt{1+\log(m)}+\dfrac{n-1}{2\sqrt{1+\log(m/n)}}.
			\end{split}
		\end{equation}
		Combining \eqref{cal01} and \eqref{cal02},
		\begin{equation}\label{cal03}
			\begin{split} 
				S_n&\le \sqrt{2}\left(n\sqrt{1+\log(m/n)}+\dfrac{n-1}{2\sqrt{1+\log(m/n)}}\right)\\
				&\le \dfrac{3n}{2}\sqrt{2(1+\log(m/n))},
			\end{split}
		\end{equation}
		which proves the second half of \eqref{lem31b}. The first half of \eqref{lem31b} coincides with \eqref{22}.
		
		Using \eqref{cal03}, we calculate lower bounds of $x_n$ and $S_n$ for $1\le n\le \lfloor m/e^3\rfloor$ as follows.
		\begin{equation}\label{cal04}
			\begin{split} 
				x_n&= \sum_{j=n}^{m-1}(x_{j}-x_{j+1})+x_m \ge\sum_{j=n}^{m-1}\dfrac{1}{S_j}\\
				&\ge \dfrac{\sqrt{2}}{3}\sum_{j=n}^{m-1} \dfrac{1}{j\sqrt{1+\log(m/j)}}\\
				&\ge \dfrac{\sqrt{2}}{3}\int_{n}^{m}\dfrac{dx}{x\sqrt{1+\log(m/x)}}\\ 
				&\ge \dfrac{1}{3}\sqrt{2(1+\log(m/n))}.
			\end{split}
		\end{equation}
		For the lower bound of $S_n$, by \eqref{cal04}, we have
		\begin{equation*}
			\begin{split} 
				S_n & \ge nx_n\ge \dfrac{n}{3}\sqrt{2\left(1+\log(m/n)\right)}.
			\end{split}
		\end{equation*}
		The proof of \eqref{lem31c} is completed.
		
		By the second half of \eqref{lem31b} and the first half of \eqref{lem31c}, we have 
		\begin{equation*}
	\begin{split} 
				& x_{i}^2-x_{j}^2=(x_{i}+x_j)(x_{i}-x_{i+1}+\dots+x_{j-1}-x_j)\\
				&= (x_{i}+x_j)\left(\dfrac{1}{S_i}+\dots+\dfrac{1}{S_{j-1}}\right)\\
				&\ge \dfrac{2\left(\sqrt{1+\log(m/{i})}+\sqrt{1+\log(m/j)}\right)}{9}\left(\sum_{k=i}^{j-1}\dfrac{1}{k\sqrt{1+\log(m/k)}}\right)\\
				&\ge \dfrac{2\left(\sqrt{1+\log(m/{i})}+\sqrt{1+\log(m/j)}\right)}{9}\int_{i}^{j}\dfrac{dx}{x\sqrt{1+\log(m/x)}}\\
				&= \dfrac{4\log(j/i)}{9}.
			\end{split}
		\end{equation*}
		This proves \eqref{lem31e}.
	\end{proof}

	\begin{proof}[Proof of \lemref{lem41}]
		We only prove \eqref{lem41_01} for the case where $x_{\lfloor m/e^3\rfloor}\le x_{n+1}< w \le x_n$. The proof for the other case is similar and therefore
		omitted.
		We have (see \cite[p. 40]{CGS})
		\begin{equation}\label{lp05}
			\begin{split} 
				g_h(w)& = \left(w-\sqrt{2\pi}(1+w^2)e^{w^2/2}(1-\Phi(w))\right)\int_{-\infty}^w h^{'}(t)\Phi(t)dt\\
				&\quad\quad -\left(w+\sqrt{2\pi}(1+w^2)e^{w^2/2}\Phi(w)\right)\int_{w}^{\infty} h^{'}(t)(1-\Phi(t))dt\\
				&:=I_1+I_2.
			\end{split}
		\end{equation}
		Applying the identity 
		\[\int_{-\infty}^w \Phi(t)dt = w\Phi(w)+\dfrac{e^{-w^2/2}}{\sqrt{2\pi}}\]
		and noting that $\|h^{'}\|=1$, we have
		\begin{equation}\label{l01}
			\begin{split} 
				\left|\int_{-\infty}^w h^{'}(t)\Phi(t)dt\right|&\le w\Phi(w)+\dfrac{e^{-w^2/2}}{\sqrt{2\pi}}\le 1+w.
			\end{split}
		\end{equation}
		Therefore, applying the first half of \eqref{lem31c} and \eqref{ap06}, and noting that $x_{\lfloor m/e^3\rfloor}\le x_{n+1}< w \le x_n$, we obtain
		\begin{equation}\label{l05}
			\begin{split} 
				|I_1|\le \dfrac{3}{(1+w)^2}\le \dfrac{C}{\log(m/n)}.
			\end{split}
		\end{equation}
		For all $1\le i\le N-1$, we have $h^{'}(t)=1$ on $(x_{i+1},m_i)$ and $h^{'}(t)=-1$ on $(m_{i},x_i)$.
		Therefore, for all $1\le i\le \lfloor m/e^3\rfloor-1$, we have
		\begin{equation}\label{lp06}
			\begin{split} 
				\left|\int\limits_{x_{i+1}}^{x_i}h^{'}(t)(1-\Phi(t))dt\right|&=\left|\int_{x_{i+1}}^{x_i}h^{'}(t)\Phi(t)dt\right|\\
				&=\left|\int_{x_{i+1}}^{m_i}\Phi(t)dt-\int_{m_{i}}^{x_i}\Phi(t)dt\right|\\
				&=\int_{m_{i}}^{x_i}\Phi(t)dt-\int_{x_{i+1}}^{m_i}\Phi(t)dt\\
				& \le\dfrac{x_i-x_{i+1}}{2}\left(\Phi(x_i)-\Phi(x_{i+1})\right)\\
				& \le \dfrac{(x_i-x_{i+1})^2 e^{-x_{i+1}^2/2}}{2\sqrt{2\pi}}\\
				&=\dfrac{e^{-x_{i+1}^2/2}}{2\sqrt{2\pi}S_{i}^2},
			\end{split}
		\end{equation}
		where we have applied the monotoncity of function $\Phi(\cdot)$ in the first inequality, the mean value theorem in the second inequality, and \eqref{MIW-recursion}
		in the last equality. 
		It follows from \eqref{lp06} that
		\begin{equation}\label{l02}
			\begin{split} 
				&\left|\int_{w}^{\infty} h^{'}(t)(1-\Phi(t))dt\right|=\left|\int_{w}^{x_1} h^{'}(t)(1-\Phi(t))dt\right|\\
				&\qquad\le \int_{w}^{x_n} (1-\Phi(t))dt+\sum_{i=1}^{n-1}\left|\int_{x_{i+1}}^{x_{i}}h^{'}(t)(1-\Phi(t)) dt\right|\\
				&\qquad\le (x_n-w)(1-\Phi(w))+\sum_{i=1}^{n-1}\dfrac{e^{-x_{i+1}^2/2}}{2\sqrt{2\pi}S_{i}^2}\\
				&\qquad\le \dfrac{(x_{n}-x_{n+1})e^{-w^2/2}}{w\sqrt{2\pi}}+\sum_{i=1}^{n-1}\dfrac{e^{-x_{i+1}^2/2}}{2\sqrt{2\pi}S_{i}^2}.
			\end{split}
		\end{equation}
		Therefore
		\begin{equation}\label{lem41_03}
			\begin{split} 
				|I_2|&= \left(w+\sqrt{2\pi}(1+w^2)e^{w^2/2}\Phi(w)\right)\left|\int_{w}^{\infty} h^{'}(t)(1-\Phi(t))dt\right|\\
				&\le C w^2e^{w^2/2}\left( \dfrac{(x_n-x_{n+1})e^{-w^2/2}}{w\sqrt{2\pi}}+\sum_{i=1}^{n-1}\dfrac{e^{-x_{i+1}^2/2}}{2\sqrt{2\pi}S_{i}^2}\right)\\
				&\le C\left(\dfrac{x_n}{S_n}+ \sum_{i=1}^{n-1}\dfrac{x_{i}^2}{S_{i}^2}e^{(x_{n}^2-x_{i+1}^2)/2}\right)\\
				&\le C\left(\dfrac{x_n}{S_n}+ \sum_{i=1}^{n-1}\dfrac{x_{i}^2}{S_{i}^2}\left(\dfrac{i+1}{n}\right)^{2/9}\right)\\
				&\le C\left(\dfrac{1}{n}+\dfrac{1}{n^{2/9}}\sum_{i=1}^{n-1}\dfrac{(i+1)^{2/9}}{i^{2}}\right)\\
				&\le \dfrac{C}{n^{2/9}},
			\end{split}
		\end{equation}
		where we have applied \eqref{l02} and the fact that $w\ge x_{\lfloor m/e^3\rfloor}\ge \sqrt{8}/3$ in the first inequality,
		\eqref{MIW-recursion} and the fact that $w^2\le x_{n}^2$ in the second inequality, \eqref{lem31e}
		in the third inequality, and the fact that $S_n\ge nx_n$ in the fourth inequality.
		The proof of \eqref{lem41_01} follows from \eqref{lp05}, \eqref{l05} and \eqref{lem41_03}.
	\end{proof}
	
\section*{Acknowledgements}
This research was partially supported by Grant  R-146-000-230-114 from the National University of Singapore. A substantial part of this paper was written when the second author was at the Vietnam Institute for Advanced Study in Mathematics (VIASM) and the Department of Mathematics, National University of Singapore (NUS). He would like to thank  VIASM and the Department of Mathematics at NUS for their hospitality. We would also like to thank Adrian R\"{o}llin for helping us draw Figure 1 and the referees for their valuable comments.

\end{document}